\documentclass{amsart}
\usepackage{amsfonts}
\usepackage{amscd}
\usepackage{amssymb}
\newtheorem{theorem}{Theorem}
\theoremstyle{plain}

\newtheorem{corollary}[theorem]{Corollary}

\theoremstyle{definition}

\newcommand{\R}{{\mathbb R}}

\newcommand{\Z}{{\mathbb Z}}
\newcommand{\C}{{\mathbb C}}

\newcommand{\D}{{\mathcal D}}

\renewcommand{\H}{{\mathcal H}}

\newcommand{\F}{{\mathcal F}}

\newcommand{\sinc}{{\mathop{\mathrm{sinc}}\,}}

\newcommand{\supp}{{\operatorname {supp}\,}}
\normalbaselines
\begin{document}

\title[Uncertainty Principles]
{{A connection between the Uncertainty Principles on the real line and on the circle}}

\author{Nils Byrial Andersen}
\address{Department of Mathematics,
Aarhus University,
Ny Munkegade 118,
Building 1530,
DK-8000 Aarhus C,
Denmark}
\email{byrial@imf.au.dk}
\date{17 July, 2013}

\begin{abstract}
The purpose of this short note is to exhibit a new connection between the Heisenberg Uncertainty Principle on the line and the
Breitenberger Uncertainty Principle on the circle, 
by considering the commutator of the multiplication and difference operators on Bernstein functions.
\end{abstract}

\maketitle

\section{Introduction}

Consider the Bernstein space $B^2_R$ defined as the subspace 
of functions $f$ in $L^2(\R)$ whose (distributive) Fourier transform $\F f$ has support in the interval $[-R,R]$. 
Let $A_\delta$ denote the family of normalized backward difference operators
\begin{equation*}
A_\delta f(z) = \frac {f(z) - f(z-\delta)}{\delta} \qquad (f \in B^2_R, \,z\in \mathbb{C}),
\end{equation*}
for $\delta\in (0,1]$. 
Let $\dot{B}^2_R= \{f\in B^2_R : x f(x) \in L^2 (\R)\}$ (note that $f\in \dot{B}^2_R \Rightarrow xf(x)\in B^2_R$), and let 
$B: \dot{B}^2_R \to B^2_R$  denote  the multiplication operator 
\begin{equation*}
B f(x) = x f(x)  \qquad (f \in \dot{B}^2_R, \,x\in \mathbb{R}).
\end{equation*}
Using an operator theoretic approach, see \cite{Erb, FS, GG1, GG2, GM, PQ, PQRS, Se}, we get the 
uncertainty inequalities
\begin{equation*}
\left\| (A_\delta -a)f\right\|_2 \left\| (B-b) f\right\| _2\ge  \frac 12 |\langle f(\cdot - \delta), f\rangle |,
\end{equation*}
for  $0\ne f\in \dot{B}^2_R$, and all $a,\,b\in \C$.

At the limit $\delta \to 0$ (and $R$ arbitrary), 
we recover the Heisenberg Uncertainty Principle ({\cite{Heis}) for functions on the line, and at
the endpoint $\delta =1$ (with $R=\pi$), we recover the 
Breitenberger Uncertainty Principle ({\cite{Br}) for functions on the circle,
thus giving a new, and easy, link between the two Uncertainty Principles. 
Another connection between the two Uncertainty Principles was 
discussed in \cite{PQRS}.

Finally, we show the equivalence of the Heisenberg Uncertainty Principle to another Uncertainty Principle on the circle
(\cite[\S3]{GG1}), using the central difference operators
\begin{equation*}
C_\delta f(z) = \frac {f(z+\delta) - f(z-\delta)}{2\delta} \qquad (f \in B^2_R, \,z\in \mathbb{C}),
\end{equation*}
for $\delta\in (0,1]$. 

We do not discuss (asymptotical) optimal functions for the Uncertainty Principles here, 
but refer to the references for a discussion of this matter. 
For more references to the subject of this article,  we refer to the recent monograph \cite{Erb}.

\section{Uncertainty Principles for symmetric and normal operators}

Let $\H$ be a Hilbert space with inner product $\langle \cdot,\cdot \rangle$ and norm $\|\cdot \| =\langle \cdot,\cdot \rangle^{1/2}$.
For $A$ and $B$ linear operators with domains $\D(A),\D(B)$ respectively, and range in $\H$, the (normalized) expectation value of the operator $A$ with respect to $f\in \D(A)$ is defined as
\begin{equation*}
\tau _A (f) := \frac {\langle A f, f\rangle}{\langle  f, f\rangle},
\end{equation*}
and the standard deviation, or variance, of the operator $A$ with respect to $f\in \D(A)$ is defined as
\begin{equation*}
\sigma _A (f) := \| Af - \tau _A (f) f\|= \min _{a\in \C} \| (A-a)f \|.
\end{equation*}
We notice that $\tau _A (f) f$ is the orthogonal projection of $A f$ onto $f$.
The commutator of $A$ and $B$ is defined as $[A,B] := AB - BA$, with domain $D(AB)\cap \D(BA)$. 

From \cite[Corollary~1]{GM}, \cite[Theorem~3.1]{Se} and \cite[Corollary~3.3]{Se}, we get the following uncertainty principle:
\begin{theorem}\label{UP}
If $A,\,B$ are symmetric or normal operators on a Hilbert space $\H$, then
\begin{equation*}
\| (A-a)f \|\| (B-b)f \|\ge \sigma _A (f) \sigma _B (f)\ge \frac 12 | \langle [A,B] f, f\rangle|,
\end{equation*}
for all nonzero $f\in  \D(AB)\cap  \D(BA)$, and all $a,\,b\in \C$.
\end{theorem}
 
Theorem~\ref{UP} can be used to prove the Heisenberg Uncertainty Principle on the real line 
(with $\H = L^2 (\R),\,A f(x) = x f(x),\,Bf(x) = i f'(x)$), 
and the Breitenberger Uncertainty Principle on the circle 
(with $\H = L^2 (-\pi,\pi),\,A f(x) = e^{ix} f(x),\,Bf(x) = i f'(x)$), see the references.

\section{The Bernstein spaces $B^2_R$}

Let $R>0$. Recall the definition of the Bernstein spaces $B^2_R$ from the introduction.
By the classical Paley--Wiener theorem, $B^2_R$ can be identified with the space of entire functions $f$ on $\C$ of exponential type 
$R$, whose restriction to $\R$ belongs to $L^2(\R)$.
It is well-known that $B^2_R$ is a Hilbert space equipped with the $L^2(\R)$ norm $\| f\|_2$ (of the restriction of $f$ to the real line),
which is invariant under differentiation, and the Bernstein inequality $\| f ' \|_2 \le R \| f\|_2$ holds, for all $f\in B^2_R$.
See also \cite[Lecture~20]{Le} for general results concerning $L^p$-Bernstein spaces.

Let $\sinc (z) =\sin(\pi z)/\pi z$.
Let $l^2(\Z)$ denote the space of square-summable sequences defined on the integers $\Z$.
Then $B^2_\pi$ and $l^2(\Z)$ are isomorphic, with proportional norms, and
the isomorphism is given by the Whittaker--Kotel'nikov--Shannon Sampling Formula
\begin{equation*}
f(z)=\sum _{n\in \Z} a_k \,{\sinc \left ({ z} -n\right )}, \qquad (\{ a_k\}_{k\in\Z}\in l^2(\Z)),
\end{equation*}
which converges in $L^2$ to a function $f\in B^2_\pi$, 
given as the unique solution of the interpolation problem $f(k) = a_k,\,k\in\Z$. 
Conversely, for any function $f \in B^2_\pi$, the sequence $\{f(k)\}_{k\in \Z}$ belongs to
$l^2(\Z)$. 

\section{Uncertainty Principles for Bernstein spaces}

Consider the family of difference operators $A_\delta: B^2_R \to B^2_R$, with $\delta\in (0,1]$, from the introduction.
For $\delta =1$, this is the usual backward difference operator $\partial f(z) = A_1 f(z) = f(z) - f(z-1)$.
The adjoint operator $A_\delta^*: B^2_R\to B^2_R$, is given by $A_\delta ^* f(z) = (f(z) - f(z+\delta))/{\delta}$, 
and since $A_\delta A_\delta^*=A_\delta^*A_\delta$, we see that $A$ is a normal operator.
The multiplication operator $B$  is obviously a symmetric operator.

A small computation yields that
\begin{equation*}
[A_\delta ,B ]f(x)= f(x-\delta)  \qquad (f \in  \dot{B}^2_R, \,x\in \mathbb{R}).
\end{equation*}
From Theorem~\ref{UP}, we thus have
\begin{theorem}\label{UPBernstein}
Let $\delta \in (0,1]$. Let $0\ne f \in \dot{B}^2_R$. Then
\begin{equation*}
\left\| (A_\delta -a)f\right\|_2 \left\| (B-b) f\right\| _2\ge   
\sigma _{A_\delta} (f) \sigma _B (f)\ge  \frac 12  \left| \left\langle f(\cdot - \delta), f\right\rangle \right|,
\end{equation*}
for all $a,\,b\in \C$.
\end{theorem}

Let $\{f(n)\}_{n\in \Z}$ be a sequence in $l^2(\Z)$, and denote by $\check{f}\in L^2(-\pi,\pi)$ the Fourier inverse of $f$. 
The Fourier series corresponding to the function $e^{i\theta}\check{f}$ is the sequence $\{f(n-1)\}$. The Fourier inverse of
$xf(x) \in {B}^2_{\pi}$, or $\{nf(n)\}\in l^2(\Z)$, is $i\frac d{d\theta} \check{f}$.
Let $\delta = 1$ and $R=\pi$, then Theorem~\ref{UPBernstein} yields
\begin{corollary}
The Breitenberger Uncertainty Principle for functions on the circle.
Let $f \in \dot{B}^2_{\pi}$. Then
\begin{equation*}
\left \| (e^{i\theta}-a)\check{f} \right \|_2 \left\| \left (\frac d{d\theta}-b\right ) \check{f} \right\|_2 
\ge \frac 12 \left |\left\langle e^{i\theta} \check{f},  \check{f}\right\rangle\right |,
\end{equation*}
for all $a,\,b\in \C$.
\end{corollary}
We can rewrite this in terms of the Fourier coefficients $\{f(n)\}$,
\begin{equation*}
\left (\sum _{n\in \Z} \left |f(n-1)-af(n)\right |^2 \right )^{\frac 12}
\left (\sum _{n\in \Z} \left | (n-b) f(n) \right |^2\right )^{\frac 12}
\ge \frac 12 \left |\sum _{n\in \Z} f(n-1)f(n) \right |,
\end{equation*}
which holds for all square-summable sequences $\{f(n)\}_{n\in \Z}\in l^2(\Z)$.

Let $\delta \to 0$, then Theorem~\ref{UPBernstein} yields
\begin{corollary}
The Heisenberg Uncertainty Principle for functions on the line.
Let $f \in {B}^2_R$, for some $R>0$. Then
\begin{equation*}
\left\|  \left (\frac d{dx} -a\right )f\right\|_2\left \| (x-b) f\right\|_2 \ge \frac 12 \left \|f\right\|_2 ^2 ,
\end{equation*}
and all $a,\,b\in \C$.
\end{corollary}
The inequality holds for $f\in \dot {B}^2_R$ by Theorem~\ref{UPBernstein}, 
and easily extends to all $f\in {B}^2_R$.
A standard density argument furthermore extends the inequality to all functions $f\in L^2(\R)$.

Finally, let us consider the central difference operators $C_\delta$. Since
\begin{equation*}
[C_\delta ,B ]f(x)= \frac {f(x+\delta)+f(x-\delta)}2  \qquad (f \in  \dot{B}^2_R, \,x\in \mathbb{R}),
\end{equation*}
Theorem~\ref{UP} gives
\begin{theorem}\label{UPBernstein2}
Let $\delta \in (0,1]$. Let  $0\ne f \in \dot{B}^2_R$. Then
\begin{equation*}
\left\| (C_\delta -a)f\right\|_2 \left\| (B-b) f\right\| _2\ge   
\sigma _{C_\delta} (f) \sigma _B (f)\ge  \frac 12 \left |\left \langle \frac {f(\cdot + \delta)+ f(\cdot - \delta)}2, f\right\rangle \right |,
\end{equation*}
for all $a,\,b\in \C$.
\end{theorem}

In the limit $\delta \to 0$, Theorem~\ref{UPBernstein2} yields
the Heisenberg Uncertainty Principle as before. So let $\delta = 1$ and $R=\pi$, 
then,
\begin{corollary}
Let $f \in \dot{B}^2_{\pi}$. Then
\begin{equation*}
\left \| (\sin (\theta)-a)\check{f} \right \|_2 \left\| \left (\frac d{d\theta}-b\right ) \check{f} \right\|_2 
\ge \frac 12 \left |\left \langle \cos (\theta) \check{f},  \check{f}\right\rangle \right |,
\end{equation*}
for all $a,\,b\in \C$.
\end{corollary}

\section{Final remarks}

Normally, when we discuss Uncertainty Principles mathematically, 
we say that $f$ or $\F f$ cannot be localized at the same time. 
Here, we assume that $\F f$ is localized as $\supp \F f \subset [-R,R]$, or $f\in {B}^2_R$, which is another reason 
why it may be interesting to look at ${B}^2_R$. 

The Bernstein inequality $\| f ' \|_2 \le R \| f\|_2$, together with the Heisenberg Uncertainty Principle, also 
yields the following inequality, for $0\ne f \in {B}^2_R$, and $a\in \C$,
\begin{equation*}
\frac 12 \left\|f\right\|_2 ^2 \le \left \|  f' \right\|_2 \left \|(x-a)  f\right \|_ 2\le R\left \| f \right\|_ 2 \left  \|(x-a) f\right\|_ 2 ,
\end{equation*}
or 
\begin{equation*}
 \left  \|(x-a) f\right\|_ 2 \ge \frac {\left\|f\right\|_2}{2R},
\end{equation*}
which indeed supports the claim that localization of the frequency, i.e.,\,$R$ small, implies indeterminacy of the position.

\end{document}